\documentclass[leqno,12pt]{article}
\usepackage{amsmath,amssymb,amsfonts,amsthm}
\usepackage{cases}
\usepackage{graphics}
\usepackage{parskip}
\usepackage{indentfirst}
\allowdisplaybreaks

\newtheorem{remark}{Remark}

\newtheorem{defin}{Definition}[section]

\title{Space-time holomorphic solutions of Navier-Stokes equations}
\author{E. Tsyganov \\
\small Department of Mathematics and Informatics Technology \\
\small Bashkir State University\\
\small entsyganov@yahoo.com} 
\date{}
\begin{document}
\maketitle
\begin{center}
\Large \bf   
\end{center}
\begin{abstract}
We introduce a concept of space-time holomorphic solutions of partial differential equations and construct a meromorphic solution of Navier-Stokes equations.
\end{abstract}
\section{Introduction}
We study Navier-Stokes equations in Lagrangean coordinates
\begin{numcases}{}
v_{t}-u_{x}=0,\label{1}\\
u_{t}+p_{x}=\left(\displaystyle \frac{u_{x}}{v}\right)_{x},\label{2} 
\end{numcases}
with Cauchy data
\begin{align}
v(0,x)=v_{0}(x), \ \ \ \ u(0,x)=u_{0}(x).
\end{align}
Here \(t\in \mathbb{R}^{+}\) and  \(x\in \mathbb{R}\) are time and space respectively, dependent variable \(v=v(x,t)\) denotes the specific volume, \(u=u(t,x)\) - velocity, \(p=p(v)\) - pressure. We assume that \(p\) satisfies the following conditions:  

$$p'<0, \ \ \ \ \lim_{v\rightarrow 0+}p=+\infty, \ \ \ \ \lim_{v\rightarrow +\infty}p=0.$$
In addition we assume that  $p$ is holomorphic  in a neigborhood of \(\mathbb{R^{+}}.\)\par
In this paper, we continue our study of solutions of the Navier-Stokes equations having analyticity properties. The issue of analyticity was first addressed in Masuda \cite{1} for Navier-Stokes equations for incompressible fluids and was further investigated in a number of papers (see Foias and Temam \cite{2}, Constantin, Foias, Kukavica, Majda \cite{3} etc.). The results show that analyticity arises naturally when solving the equations in their classical form, and can be used to study their properties.\\
\indent Of special interest is the study of complex solutions of the Navier-Stokes equations. Despite the fact that these models do not have a direct physical significance, they provide new information about the equations themselves (see Dong Li and Y. Sinai \cite{4}).\\
\indent The study of analytic properties of weak solutions of the Navier-Stokes equations for a compressible gas dynamic was initiated in E.Tsyganov \cite{5} and was further developed for the multi-dimensional case in D.Hoff and E.Tsyganov \cite{6}. We can point out that analyticity plays a critical role in the proof of backward uniqueness and in the derivation of exact rates of regularization and asymptotic behavior of weak solutions (see E.Tsyganov \cite{7}).\\
\indent In this paper, we study a special class of analytic solutions, which we call space-time holomorphic. The basic idea is to merge $t$ and $ x $ in one complex variable by the equality $z = it + x$, where $ i $ is the imaginary unit. After that, in order to find solutions of the original partial differential equation, we need to solve an ordinary differential equation in the complex plane. We give an example of a family of solutions, which are  elementary meromorphic functions on the whole complex plane. These functions may blow-up (have a pole) in  finite time, however, they become smooth again once we  go through these positive values of time. We also prove that these solutions are space and time meromorphic. \\
\indent The article is structured as follows: in Section 2  we introduce the notion of space-time holomorphic solutions of partial differential equations. Then, in Section 3, we give an example of such a solution of the Navier-Stokes equations and study its properties.
\section{Space-time holomorphic solutions}
In this section we introduce a concept of space-time holomorphic solutions.
\begin{defin}
 Consider a partial differential equation 
\begin{align}
F\left(u,\frac{\partial u}{\partial t},...,\frac{\partial^{n} u}{\partial t^{n}},\frac{\partial u}{\partial x},...,\frac{\partial^{m} u}{\partial x^{m}}\right)=0,\label{4}
\end{align}
where $t\in \mathbb{R}, x\in \mathbb{R}, u\in \mathbb{R}^{l}$, and $F\in \mathbb{R}^{k}$ is a holomorphic function of its arguments.\par
We say that solution $u$ of equation \eqref{4} is space-time holomorphic, if, in addition, it satisfies the equation  
\begin{align}
\frac{\partial u}{\partial t}=i\frac{\partial u}{\partial x},\label{5}
\end{align} 
where $i$ is the imaginary unit. 
\end{defin}
Then it follows from condition \eqref{5} that function $u$ satisfies Cauchy-Riemann condition. So, if we set $z=x+it$, then $u$ becomes a holomorphic function of $z$ and the following equalities hold:
\begin{align}
\frac{\partial u}{\partial t}=i\frac{d u}{d z}, \ \ \ \ \ \  \frac{\partial u}{\partial x}=\frac{d u}{d z}.
\end{align} 
\indent Now we can give another definition of space-time holomorphy.
\begin{defin}
 We say that a function $u$ is a space-time holomorphic solution of equation \eqref{4}, if it satisfies an ordinary differential equation 
\begin{align}
F\left(u,i\frac{d u}{d z},...,i^{n}\frac{d^{n} u}{d z^{n}},\frac{d u}{d z},...,\frac{d^{m} u}{d z^{m}}\right)=0\label{6}
\end{align} 
in some domain of the complex plane.
\end{defin}
\begin{remark}
Instead of condition \eqref{5} we can set
\begin{align}
\frac{\partial u}{\partial t}=i\frac{\partial u}{\partial x}\tag{4*}\label{13}.
\end{align} 
If we now let $z=t+ix$, then we will come to a different equation for the space-time holomorphy:
\begin{align}
F\left(u,\frac{d u}{d z},...,\frac{d^{n} u}{d z^{n}},i\frac{d u}{d z},...,i^{m}\frac{d^{m} u}{d z^{m}}\right)=0 \tag{\ref{6}*}.
\end{align} 
\end{remark}
Throughout the rest of the paper we will only consider equation \eqref{6}.
\begin{remark}
Let $u=u(z)$ satisfy the equation \eqref{6} in some domain of the complex plane. If we set $x=Re(z), \ t=Im(z)$, then $u=u(x,t)$ satisfies partial differential equation \eqref{4} in some domain in the plane $OXY$.
\end{remark}
\section{Meromorphic solutions of Navier-Stokes equations}
The equations of the space-time holomorphy \eqref{6} for Navier-Stokes equations \eqref{1},\eqref{2} are the following:
\begin{numcases}{}
iv_{z}-u_{z}=0;\\
iu_{z}+p_{z}=\left(\frac{u_{z}}{v}\right)_{z}.
\end{numcases}
We express $v_{z}$ in terms of $u_{z}$ and then substitute it into the second equation:
\begin{align}
-v_{z}+p_{z}=\left(\frac{iv_{z}}{v}\right)_{z}.
\end{align}
Then we integrate the last equality to obtain 
\begin{align}
-v+p=i\frac{v_{z}}{v}+C,\label{13}
\end{align}
where $C$ is an arbitrary constant of integration. We rewrite the equation and integrate it once again:
\begin{align}
z+C_{1}=i\int \frac{dv}{-v^{2}+vp+C_{2}v}
\end{align}
We will primarily be interested in those functions $ p $, for which the equation is integrable by quadratures.
Setting $p=\dfrac{1}{v}$, $C_{2}=0$ and $C_{1}=C$, we can obtained the answer in closed form:
$$z+C=i\frac{1}{2}(\ln(v+1)-\ln(1-v)),$$
from which we obtain $v$ :
\begin{align}
v=\frac{e^{-2i(z+C)}-1}{e^{-2i(z+C)}+1}.
\end{align}
Now we can find $u$ :
\begin{align}
u=i\frac{e^{-2i(z+C)}-1}{e^{-2i(z+C)}+1}+C_{2}.
\end{align}
We point out that the functions $v$ and $u$ are meromorphic  on the whole complex plane and have the following properties:
\begin{align}
\lim_{\stackrel{z\rightarrow \infty,}{ Im(z)>0}} v=1, \ \ \lim_{\stackrel{z\rightarrow \infty,}{ Im(z)>0}} u=i+C_{2}.
\end{align}
\indent We are interested primarily in the zeros and poles of $v$:
\begin{align*}
\text{zeroes : } z=\pi k-C, \ \ \text{poles : } z=\frac{\pi}{2}+\pi k -C, \ \ k\in \mathbb{Z}.
\end{align*} 
\indent Depending on the value of $ C $ we can  have the following situation in the half-plane $ Im (z) \geq 0 $:\newline
1) $Im(C)>0$. Then the functions $ v $ and $ u $ are holomorphic in $ Im (z)> 0 $ and continuous on $ Im (z) \geq 0 $.
In this case, the pair of functions $ v = v (x, t) $, $ u = u (x, t) $ is a smooth solution of the Navier-Stokes equations \eqref {1}, \eqref {2} in the half-plane $ \mathbb {R} \times \mathbb {R} ^ {+} $ with smooth initial data; \newline
2) $Im(C)=0$. Then $ (v, u) $ is a smooth solution of \eqref {1},\eqref {2} in the half-plane $t> 0$ with meromorphic initial data;\newline
3) $Im(C)<0$. In this case, $ (v, u) $ is a meromorphic solution for $ t> 0 $ with smooth initial data at $ t = 0 $. 
\begin{remark}
Equation \eqref{2} will have singularity at the points where $v=0$. These singularities are, however, removable, since $(v,u)$ exist and is smooth at  such points.
\end{remark}
\begin{remark}
\indent Functions $ v = v (t, x), u = u (t, x) $ are meromorphic in  $ t $ for any fixed $ x \ in \mathbb{R} $.  To prove this we use the following identities:
\begin{align}
v(t,x)\stackrel{\, def}{=\joinrel=}v(t+x), \ \ \ \ u(t,x)\stackrel{\, def}{=\joinrel=}u(t+x), \ \ \ \ t\in \mathbb{C}, \ x\in \mathbb{R}. 
\end{align} 
The solution is also meromorphic in $ x $ for any fixed $ t \in \mathbb{R} $. This follows from the obvious definitions
\begin{align}
v(t,x)\stackrel{\, def}{=\joinrel=}v(it+x), \ \ \ \ u(t,x)\stackrel{\, def}{=\joinrel=}u(it+x), \ \ \ \ x\in \mathbb{C}, \ t\in \mathbb{R}.
\end{align}
\end{remark}

\end{document}